# Similar Ruled Surfaces with Variable Transformations in the Euclidean 3-space $E^3$

**Mehmet Önder**
*Celal Bayar University, Faculty of Arts and Sciences, Department of Mathematics, Muradiye Campus, 45047 Muradiye, Manisa, Turkey. E-mail*: mehmet.onder@cbu.edu.tr

**Abstract**
In this study, we define a family of ruled surfaces in the Euclidean 3-space $E^3$ and called similar ruled surfaces. We obtain some properties of these special surfaces and we show that developable ruled surfaces form a family of similar ruled surfaces if and only if the striction curves of the surfaces are similar curves with variable transformation.

**MSC:** 53A25, 14J26.

**Key words:** Ruled surface; Frenet frame; similar curves.

## 1. Introduction

In local differential geometry, associated curves, the curves for which at the corresponding points of the curves one of the Frenet vectors of a curve coincides with one of the Frenet vectors of other curve, are very interesting and an important problem of the fundamental theory and the characterizations of space curves. The well-known pairs of such curves are Bertrand curves, Mannheim curves and involute-evolute curves [4,9,10]. Recently, a new definition of the associated curves was given by El-Sabbagh and Ali [1]. They have called these new curves as similar curves with variable transformation and defined as follows: Let $\psi_\alpha(s_\alpha)$ and $\psi_\beta(s_\beta)$ be two regular curves in $E^3$ parametrized by arc lengths $s_\alpha$ and $s_\beta$ with curvatures $\kappa_\alpha$, $\kappa_\beta$ and torsions $\tau_\alpha$, $\tau_\beta$ and Frenet frames $\{\vec{T}_\alpha, \vec{N}_\alpha, \vec{B}_\alpha\}$ and $\{\vec{T}_\beta, \vec{N}_\beta, \vec{B}_\beta\}$, respectively. $\psi_\alpha(s_\alpha)$ and $\psi_\beta(s_\beta)$ are called similar curves with variable transformation $\lambda_\beta^\alpha$ if there exists a variable transformation

$$s_\alpha = \int \lambda_\beta^\alpha(s_\beta) ds_\beta,$$

of the arc lengths such that the tangent vectors are the same for two curves i.e., $\vec{T}_\alpha = \vec{T}_\beta$ for all corresponding values of parameters under the transformation $\lambda_\beta^\alpha$. All curves satisfying this condition is called a family of similar curves. Moreover, they have obtained some properties of the family of similar curves.

Analogue to the special curve pairs, the surface pairs, especially ruled surface pairs (called offset surfaces), have an important positions and applications in the study of design problems in spatial mechanisms and physics, kinematics and computer aided design (CAD) [6,7]. So, these surfaces are one of the most important topics of surface theory. In fact, ruled surface offsets are the generalization of the notion of Bertrand curves and Mannheim curves to line geometry and these surface pairs are called Bertrand offsets and Mannheim offsets [3,5,8].

In this work, we consider the notion of similar curves for ruled surfaces. We introduce a family of ruled surfaces in the Euclidean 3-space $E^3$ and called similar ruled surfaces with variable transformation. We give some theorems characterizing these special surfaces and we show that developable ruled surfaces form a family of similar ruled surfaces if and only if the striction curves of the surfaces are similar curves with variable transformation.

## 2. Ruled Surfaces in $E^3$

In this section, we give a brief summary of the theory of ruled surface in $E^3$. A more detailed information can be obtained in ref. [2].

Let $I$ be an open interval in the real line $IR$, $\vec{f} = \vec{f}(u)$ be a curve in $E^3$ defined on $I$ and $\vec{q} = \vec{q}(u)$ be a unit direction vector of an oriented line in $E^3$. Then we have the following parametrization for a ruled surface $N$,

$$\vec{r}(u,v) = \vec{f}(u) + v\vec{q}(u). \tag{1}$$

The curve $\vec{f} = \vec{f}(u)$ is called base curve or generating curve of the surface and various positions of the generating lines $\vec{q} = \vec{q}(u)$ are called rulings. In particular, if the direction of $\vec{q}$ is constant, then the ruled surface is said to be cylindrical, and non-cylindrical otherwise.

The distribution parameter of $N$ is given by

$$d = \frac{\left|\dot{\vec{f}}, \vec{q}, \dot{\vec{q}}\right|}{\left\langle \dot{\vec{q}}, \dot{\vec{q}} \right\rangle}, \tag{2}$$

where $\dot{\vec{f}} = \frac{d\vec{f}}{du}$, $\dot{\vec{q}} = \frac{d\vec{q}}{du}$. If $\left|\dot{\vec{f}}, \vec{q}, \dot{\vec{q}}\right| = 0$, then normal vectors are collinear at all points of same ruling and at nonsingular points of the surface $N$, the tangent planes are identical. We then say that tangent plane contacts the surface along a ruling. Such a ruling is called a torsal ruling. If $\left|\dot{\vec{f}}, \vec{q}, \dot{\vec{q}}\right| \neq 0$, then the tangent planes of the surface $N$ are distinct at all points of same ruling which is called nontorsal.

**Definition 2.1.** ([2]) A ruled surface whose all rulings are torsal is called a developable ruled surface. The remaining ruled surfaces are called skew ruled surfaces. From (2) it is clear that a ruled surface is developable if and only if at all its points the distribution parameter is zero.

For the unit normal vector $\vec{m}$ of a ruled surface $N$ we have

$$\vec{m} = \frac{\vec{r}_u \times \vec{r}_v}{\|\vec{r}_u \times \vec{r}_v\|} = \frac{(\dot{\vec{f}} + v\dot{\vec{q}}) \times \vec{q}}{\sqrt{\left\langle \dot{\vec{f}} + v\dot{\vec{q}}, \dot{\vec{f}} + v\dot{\vec{q}} \right\rangle - \left\langle \dot{\vec{f}}, \vec{q} \right\rangle^2}}. \tag{3}$$

The unit normal of the surface along a ruling $u = u_1$ approaches a limiting direction as $v$ infinitely decreases. This direction is called the asymptotic normal (central tangent) direction and from (3) defined by

$$\vec{a} = \lim_{v \to \pm\infty} \vec{m}(u_1, v) = \frac{\vec{q} \times \dot{\vec{q}}}{\|\dot{\vec{q}}\|}.$$

The point at which the unit normal of $N$ is perpendicular to $\vec{a}$ is called the striction point (or central point) $C$ and the set of striction points on all rulings is called striction curve of the surface. The parametrization of the striction curve $\vec{c} = \vec{c}(u)$ on a ruled surface is given by

$$\vec{c}(u) = \vec{f}(u) + v_0\vec{q}(u) = \vec{f} - \frac{\left\langle \dot{\vec{q}}, \dot{\vec{f}} \right\rangle}{\left\langle \dot{\vec{q}}, \dot{\vec{q}} \right\rangle} \vec{q}, \tag{4}$$

where $v_0 = -\dfrac{\langle \dot{\vec{q}}, \dot{\vec{f}} \rangle}{\langle \dot{\vec{q}}, \dot{\vec{q}} \rangle}$ is called strictional distance.

The vector $\vec{h}$ defined by $\vec{h} = \vec{a} \times \vec{q}$ is called central normal which is the surface normal along the striction curve. Then the orthonormal system $\{C; \vec{q}, \vec{h}, \vec{a}\}$ is called Frenet frame of the ruled surface $N$ where $C$ is the central point of ruling of ruled surface $N$ and $\vec{q}$, $\vec{h} = \vec{a} \times \vec{q}$, $\vec{a}$ are unit vectors of ruling, central normal and central tangent, respectively.

For the derivatives of the vectors of Frenet frame $\{C; \vec{q}, \vec{h}, \vec{a}\}$ of ruled surface $N$ with respect to the arc length $s$ of striction curve we have

$$\begin{bmatrix} d\vec{q}/ds \\ d\vec{h}/ds \\ d\vec{a}/ds \end{bmatrix} = \begin{bmatrix} 0 & k_1 & 0 \\ -k_1 & 0 & k_2 \\ 0 & -k_2 & 0 \end{bmatrix} \begin{bmatrix} \vec{q} \\ \vec{h} \\ \vec{a} \end{bmatrix} \qquad (5)$$

where $k_1 = \dfrac{ds_1}{ds}$, $k_2 = \dfrac{ds_3}{ds}$ and $s_1$, $s_3$ are the arc lengths of the spherical curves circumscribed by the bound vectors $\vec{q}$ and $\vec{a}$, respectively. The ruled surfaces satisfying $k_1 \neq 0$, $k_2 = 0$ are called conoids (For details [2]).

Now, we can represent and prove the following theorems which are necessary for the following section.

**Theorem 2.1.** *Let the striction curve $\vec{c} = \vec{c}(s)$ of ruled surface $N$ be unit speed i.e., $s$ is arc length parameter of $\vec{c}(s)$ and let $\vec{c}(s)$ be the base curve of the surface. Then $N$ is developable if and only if the unit tangent of the striction curve is the same with the ruling along the curve.*

**Proof:** Let $s$ be arc length parameter of the striction curve. Then the unit tangent of the striction curve is given by

$$\vec{T}(s) = \frac{d\vec{c}}{ds} = (\cos\theta)\vec{q}(s) + (\sin\theta)\vec{a}(s),$$

where $\theta = \theta(s)$ is the angle between unit vectors $\vec{T}(s)$ and $\vec{q}(s)$. Since the striction curve is base curve, then from (2) and (5) the distribution parameter of the surface $N$ is obtained as

$$d = \frac{\sin\theta}{k_1}.$$

Thus we have that $N$ is developable if and only if $\sin\theta = 0$, i.e., $\vec{T}(s) = \vec{q}(s)$ satisfies.

**Theorem 2.2.** *Let the striction curve $\vec{c} = \vec{c}(s)$ of ruled surface $N$ be unit speed i.e., $s$ is arc length parameter of $\vec{c}(s)$. Suppose that $\vec{c} = \vec{c}(\varphi)$ is another parametrization of striction curve by the parameter $\varphi(s) = \int k_1(s)ds$. Then the ruling $\vec{q}$ satisfies a vector differential equation of third order given by*

$$\frac{d}{d\varphi}\left(\frac{1}{f(\varphi)} \frac{d^2 \vec{q}}{d\varphi^2}\right) + \left(\frac{1 + f^2(\varphi)}{f(\varphi)}\right) \frac{d\vec{q}}{d\varphi} - \left(\frac{1}{f^2(\varphi)} \frac{df(\varphi)}{d\varphi}\right) \vec{q} = 0, \qquad (6)$$

where $f(\varphi) = \dfrac{k_2(\varphi)}{k_1(\varphi)}$.

**Proof:** If we write derivatives given in (5) according to $\varphi$, we have

$$\frac{d\vec{q}}{d\varphi} = \frac{d\vec{q}}{ds}\frac{ds}{d\varphi} = (k_1\vec{h})\frac{1}{k_1} = \vec{h},$$

$$\frac{d\vec{h}}{d\varphi} = \frac{d\vec{h}}{ds}\frac{ds}{d\varphi} = (-k_1\vec{q} + k_2\vec{a})\frac{1}{k_1} = -\vec{q} + f(\varphi)\vec{a},$$

$$\frac{d\vec{a}}{d\varphi} = \frac{d\vec{a}}{ds}\frac{ds}{d\varphi} = (-k_2\vec{h})\frac{1}{k_1} = -f(\varphi)\vec{h},$$

respectively, where $f(\varphi) = \dfrac{k_2(\varphi)}{k_1(\varphi)}$. Then corresponding matrix form of (5) can be given

$$\begin{bmatrix} d\vec{q}/d\varphi \\ d\vec{h}/d\varphi \\ d\vec{a}/d\varphi \end{bmatrix} = \begin{bmatrix} 0 & 1 & 0 \\ -1 & 0 & f(\varphi) \\ 0 & -f(\varphi) & 0 \end{bmatrix} \begin{bmatrix} \vec{q} \\ \vec{h} \\ \vec{a} \end{bmatrix}. \tag{7}$$

From the first and second equations of new Frenet derivatives (7) we have

$$\vec{a} = \frac{1}{f(\varphi)}\left(\frac{d^2\vec{q}}{d\varphi^2} + \vec{q}\right). \tag{8}$$

Substituting the above equation in the last equation of (7) we have desired equation (6).

**3. Similar Ruled Surfaces with Variable Transformations**

In this section we introduce the definition and characterizations of similar ruled surfaces with variable transformation in $E^3$. First, we give the following definition.

**Definition 3.1.** Let $N_\alpha$ and $N_\beta$ be two ruled surfaces in $E^3$ given by the parametrizations

$$\begin{cases} \vec{r}_\alpha(s_\alpha, v) = \vec{\alpha}(s_\alpha) + v\vec{q}_\alpha(s_\alpha), & \|\vec{q}_\alpha(s_\alpha)\| = 1, \\ \vec{r}_\beta(s_\beta, v) = \vec{\beta}(s_\beta) + v\vec{q}_\beta(s_\beta), & \|\vec{q}_\beta(s_\beta)\| = 1, \end{cases} \tag{9}$$

respectively, where $\vec{\alpha}(s_\alpha)$ and $\vec{\beta}(s_\beta)$ are striction curves of $N_\alpha$ and $N_\beta$ and $s_\alpha$, $s_\beta$ are arc length parameters of $\vec{\alpha}(s_\alpha)$ and $\vec{\beta}(s_\beta)$, respectively. Let the Frenet frames and invariants of $N_\alpha$ and $N_\beta$ be $\{\vec{q}_\alpha, \vec{h}_\alpha, \vec{a}_\alpha\}$, $k_1^\alpha, k_2^\alpha$ and $\{\vec{q}_\beta, \vec{h}_\beta, \vec{a}_\beta\}$, $k_1^\beta, k_2^\beta$, respectively. $N_\alpha$ and $N_\beta$ are called similar ruled surfaces with variable transformation $\lambda_\beta^\alpha$ if there exists a variable transformation

$$s_\alpha = \int \lambda_\beta^\alpha(s_\beta) ds_\beta, \tag{10}$$

of the arc lengths of striction curves such that the rulings are the same for two ruled surfaces i.e.,

$$\vec{q}_\alpha(s_\alpha) = \vec{q}_\beta(s_\beta), \tag{11}$$

for all corresponding values of parameters under the transformation $\lambda_\beta^\alpha$. All ruled surfaces satisfying equation (11) are called a family of similar ruled surfaces with variable transformation.

Then we can give the following theorems characterizing similar ruled surfaces. Whenever we talk about $N_\alpha$ and $N_\beta$ we mean that surfaces have the parametrizations as given in (9).

**Theorem 3.1.** *Let $N_\alpha$ and $N_\beta$ be two ruled surfaces in $E^3$. Then $N_\alpha$ and $N_\beta$ are similar ruled surfaces with variable transformation if and only if the central normal vectors of the surfaces are the same, i.e.,*

$$\vec{h}_\alpha(s_\alpha) = \vec{h}_\beta(s_\beta), \tag{12}$$

*under the particular variable transformation*

$$\lambda_\beta^\alpha = \frac{ds_\alpha}{ds_\beta} = \frac{k_1^\beta}{k_1^\alpha}, \tag{13}$$

*of the arc lengths.*

**Proof:** Let $N_\alpha$ and $N_\beta$ be two similar ruled surfaces in $E^3$ with variable transformation. Then differentiating (11) with respect to $s_\beta$ it follows

$$k_1^\alpha \lambda_\beta^\alpha \vec{h}_\alpha = k_1^\beta \vec{h}_\beta. \tag{14}$$

From (14) we obtain (12) and (13) immediately.

Conversely, let $N_\alpha$ and $N_\beta$ be two ruled surfaces in $E^3$ satisfying (12) and (13). By multiplying (12) with $k_1^\beta$ and differentiating the result equality with respect to $s_\beta$ we have

$$\int k_1^\beta(s_\beta) \vec{h}_\beta(s_\beta) ds_\beta = \int k_1^\beta(s_\beta) \vec{h}_\beta(s_\beta) \frac{ds_\beta}{ds_\alpha} ds_\alpha. \tag{15}$$

From (12) and (13) we obtain

$$\vec{q}_\beta(s_\beta) = \int k_1^\beta(s_\beta) \vec{h}_\beta(s_\beta) ds_\beta = \int k_1^\alpha(s_\alpha) \vec{h}_\alpha(s_\alpha) ds_\alpha = \vec{q}_\alpha(s_\alpha), \tag{16}$$

which means that $N_\alpha$ and $N_\beta$ are similar ruled surfaces with variable transformation.

**Theorem 3.2.** *Let $N_\alpha$ and $N_\beta$ be two ruled surfaces in $E^3$. Then $N_\alpha$ and $N_\beta$ are similar ruled surfaces with variable transformation if and only if the asymptotic normal vectors of the surfaces are the same, i.e.,*

$$\vec{a}_\alpha(s_\alpha) = \vec{a}_\beta(s_\beta), \tag{17}$$

*under the particular variable transformation*

$$\lambda_\beta^\alpha = \frac{ds_\alpha}{ds_\beta} = \frac{k_2^\beta}{k_2^\alpha}, \tag{18}$$

*of the arc lengths.*

**Proof:** Let $N_\alpha$ and $N_\beta$ be two similar ruled surfaces in $E^3$ with variable transformation. Then from Definition 3.1 and Theorem 3.1 there exists a variable transformation of the arc lengths such that the rulings and central normal vectors are the same. Then from (11) and (12) we have

$$\vec{a}_\alpha(s_\alpha) = \vec{q}_\alpha(s_\alpha) \times \vec{h}_\alpha(s_\alpha) = \vec{q}_\beta(s_\beta) \times \vec{h}_\beta(s_\beta) = \vec{a}_\beta(s_\beta). \tag{19}$$

Conversely, let $N_\alpha$ and $N_\beta$ be two ruled surfaces in $E^3$ satisfying (17) and (18). By differentiating (17) with respect to $s_\beta$ we obtain

$$k_2^\alpha(s_\alpha)\vec{h}_\alpha(s_\alpha)\frac{ds_\alpha}{ds_\beta} = k_2^\beta(s_\beta)\vec{h}_\beta(s_\beta), \tag{20}$$

which gives us

$$\lambda_\beta^\alpha = \frac{k_2^\beta}{k_2^\alpha}, \quad \vec{h}_\alpha(s_\alpha) = \vec{h}_\beta(s_\beta). \tag{21}$$

Then from (17) and (21) we have

$$\vec{q}_\alpha(s_\alpha) = \vec{h}_\alpha(s_\alpha) \times \vec{a}_\alpha(s_\alpha) = \vec{h}_\beta(s_\beta) \times \vec{a}_\beta(s_\beta) = \vec{q}_\beta(s_\beta), \tag{22}$$

which completes the proof.

**Theorem 3.3.** *Let $N_\alpha$ and $N_\beta$ be two ruled surfaces in $E^3$ with non-zero curvatures. Then $N_\alpha$ and $N_\beta$ are similar ruled surfaces with variable transformation if and only if the ratio of curvatures are the same i.e.,*

$$\frac{k_2^\beta(s_\beta)}{k_1^\beta(s_\beta)} = \frac{k_2^\alpha(s_\alpha)}{k_1^\alpha(s_\alpha)}, \tag{23}$$

*under the particular variable transformation keeping equal total curvatures, i.e.,*

$$\varphi_\beta(s_\beta) = \int k_1^\beta(s_\beta)ds_\beta = \int k_1^\alpha(s_\alpha)ds_\alpha = \varphi_\alpha(s_\alpha) \tag{24}$$

*of the arc lengths.*

**Proof:** Let $N_\alpha$ and $N_\beta$ be two regular similar ruled surfaces in $E^3$ with variable transformation. Then from (18) and (21) we have (23) under the variable transformation (24), and this transformation is also leads from (18) by integration.

Conversely, let $N_\alpha$ and $N_\beta$ be two ruled surfaces in $E^3$ satisfying (23) and (24). From Theorem 2.2, the rulings $\vec{q}_\alpha$ and $\vec{q}_\beta$ of the surfaces $N_\alpha$ and $N_\beta$ satisfy the following vector differential equations of third order

$$\frac{d}{d\varphi_\alpha}\left(\frac{1}{f_\alpha(\varphi_\alpha)}\frac{d^2\vec{q}_\alpha}{d\varphi_\alpha^2}\right) + \left(\frac{1+f_\alpha^2(\varphi_\alpha)}{f_\alpha(\varphi_\alpha)}\right)\frac{d\vec{q}_\alpha}{d\varphi_\alpha} - \left(\frac{1}{f_\alpha^2(\varphi_\alpha)}\frac{df_\alpha(\varphi_\alpha)}{d\varphi_\alpha}\right)\vec{q}_\alpha = 0, \tag{25}$$

$$\frac{d}{d\varphi_\beta}\left(\frac{1}{f_\beta(\varphi_\beta)}\frac{d^2\vec{q}_\beta}{d\varphi_\beta^2}\right) + \left(\frac{1+f_\beta^2(\varphi_\beta)}{f_\beta(\varphi_\beta)}\right)\frac{d\vec{q}_\beta}{d\varphi_\beta} - \left(\frac{1}{f_\beta^2(\varphi_\beta)}\frac{df_\beta(\varphi_\beta)}{d\varphi_\beta}\right)\vec{q}_\beta = 0, \tag{26}$$

where $f_\alpha(\varphi_\alpha) = \frac{k_2^\alpha(\varphi_\alpha)}{k_1^\alpha(\varphi_\alpha)}$, $f_\beta(\varphi_\beta) = \frac{k_2^\beta(\varphi_\beta)}{k_1^\beta(\varphi_\beta)}$, $\varphi_\alpha(s_\alpha) = \int k_1^\alpha(s_\alpha)ds_\alpha$, $\varphi_\beta(s_\beta) = \int k_1^\beta(s_\beta)ds_\beta$.

From (23) we have $f_\alpha(\varphi_\alpha) = f_\beta(\varphi_\beta)$ under the variable transformation $\varphi_\alpha = \varphi_\beta$. Thus under the equation (23) and transformation (24), the equations (25) and (26) are the same, i.e., they have the same solutions. It means that the rulings $\vec{q}_\alpha$ and $\vec{q}_\beta$ are the same. Then $N_\alpha$ and $N_\beta$ are two similar ruled surfaces in $E^3$ with variable transformation.

**Theorem 3.4.** *Let the ruled surfaces $N_\alpha$ and $N_\beta$ be developable. Then $N_\alpha$ and $N_\beta$ are similar ruled surfaces with variable transformation if and only if the striction curves of the surfaces are similar curves with variable transformation.*

**Proof:** Let developable ruled surfaces $N_\alpha$ and $N_\beta$ be two similar ruled surfaces in $E^3$ with variable transformation. Since the surfaces are developable, from Theorem 2.1 we have

$$\frac{d\vec{\alpha}}{ds_\alpha} = \vec{T}_\alpha(s_\alpha) = \vec{q}_\alpha(s_\alpha), \quad \frac{d\vec{\beta}}{ds_\beta} = \vec{T}_\beta(s_\beta) = \vec{q}_\beta(s_\beta). \tag{27}$$

where $\vec{T}_\alpha(s_\alpha)$ and $\vec{T}_\beta(s_\beta)$ are unit tangents of the curves $\vec{\alpha}(s_\alpha)$ and $\vec{\beta}(s_\beta)$, respectively. From (11) and (27) we have

$$\frac{d\vec{\alpha}}{ds_\alpha} = \vec{q}_\alpha(s_\alpha) = \vec{q}_\beta(s_\beta) = \frac{d\vec{\beta}}{ds_\beta} \tag{28}$$

which shows that striction curves $\vec{\alpha}(s_\alpha)$ and $\vec{\beta}(s_\beta)$ are similar curves.

Conversely, if the striction curves $\vec{\alpha}(s_\alpha)$ and $\vec{\beta}(s_\beta)$ are similar curves, then there exists a variable transformation between arc lengths such that

$$\frac{d\vec{\alpha}}{ds_\alpha} = \vec{T}_\alpha(s_\alpha) = \vec{T}_\beta(s_\beta) = \frac{d\vec{\beta}}{ds_\beta}. \tag{29}$$

Since the ruled surfaces are developable, from Theorem 2.1 we have $\vec{T}_\alpha(s_\alpha) = \vec{q}_\alpha(s_\alpha)$ and $\vec{T}_\beta(s_\beta) = \vec{q}_\beta(s_\beta)$. Then from (29) we have that $\vec{q}_\alpha(s_\alpha) = \vec{q}_\beta(s_\beta)$, i.e., $N_\alpha$ and $N_\beta$ are similar ruled surfaces with variable transformation.

Let now consider some special cases. From (13) and (18) we have
$$k_1^\beta = \lambda_\beta^\alpha k_1^\alpha, \quad k_2^\beta = \lambda_\beta^\alpha k_2^\alpha, \tag{30}$$

respectively. From (30) it is clear that if $N_\alpha$ is a cylindrical surface i.e., $k_1^\alpha = 0$, then under the variable transformation the curvature does not change. So we have the following corollaries.

***Corollary 3.1.*** *The family of cylindrical surfaces forms a family of similar ruled surfaces with variable transformation.*

If $N_\alpha$ is a conoid surface i.e., $k_2^\alpha = 0$, then under the variable transformation the curvature does not change. So we have the following corollary.

***Corollary 3.2.*** *The family of conoid surfaces forms a family of similar ruled surfaces with variable transformation.*

**4. Example:** Let consider the ruled surface $N_\beta$ given by the parametrization

$$\begin{aligned}\vec{r}_\beta(s_\beta, v) &= \vec{\beta}(s_\beta) + v\vec{q}_\beta(s_\beta) \\ &= (0, 0, s_\beta) + v(\cos s_\beta, \sin s_\beta, 0)\end{aligned}$$

which is plotted in Fig. 1. The Frenet vectors of $N_\beta$ are

$$\vec{q}_\beta = (\cos s_\beta, \sin s_\beta, 0),$$
$$\vec{h}_\beta = (-\sin s_\beta, \cos s_\beta, 0),$$
$$\vec{a}_\beta = (0, 0, 1),$$

and curvatures are obtained as $k_1^\beta = 1$, $k_2^\beta = 0$. A similar ruled surfaces of $N_\beta$ is the surface $N_\alpha$ given by the parametrization

$$\vec{r}_\alpha(s_\alpha, v) = \vec{\alpha}(s_\alpha) + v\vec{q}_\alpha(s_\alpha)$$
$$= \left(-\sin\frac{s_\alpha}{\sqrt{2}}, \cos\frac{s_\alpha}{\sqrt{2}}, \frac{s_\alpha}{\sqrt{2}}\right) + v\left(\cos\frac{s_\alpha}{\sqrt{2}}, \sin\frac{s_\alpha}{\sqrt{2}}, 0\right)$$

which is plotted in Fig. 2. The Frenet vectors of $N_\alpha$ are

$$\vec{q}_\alpha = \left(\cos\frac{s_\alpha}{\sqrt{2}}, \sin\frac{s_\alpha}{\sqrt{2}}, 0\right),$$
$$\vec{h}_\alpha = \left(-\sin\frac{s_\alpha}{\sqrt{2}}, \cos\frac{s_\alpha}{\sqrt{2}}, 0\right),$$
$$\vec{a}_\alpha = (0, 0, 1),$$

and curvatures are obtained as $k_1^\alpha = \frac{1}{\sqrt{2}}$, $k_2^\alpha = 0$. From Theorem 3.1 we have the particular variable transformation

$$\lambda_\beta^\alpha = \frac{ds_\alpha}{ds_\beta} = \frac{k_1^\beta}{k_1^\alpha} = \sqrt{2},$$

which means that $s_\alpha = \sqrt{2} s_\beta$.

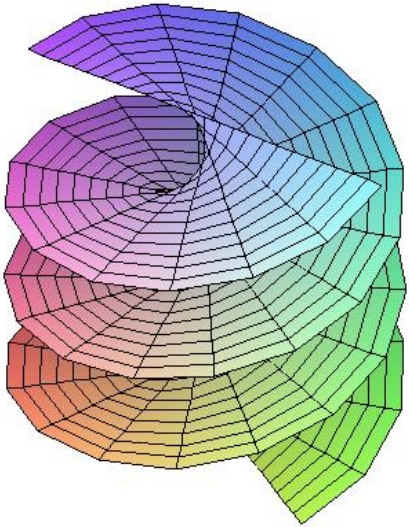
**Fig. 1.** Helicoid surface $N_\beta$.

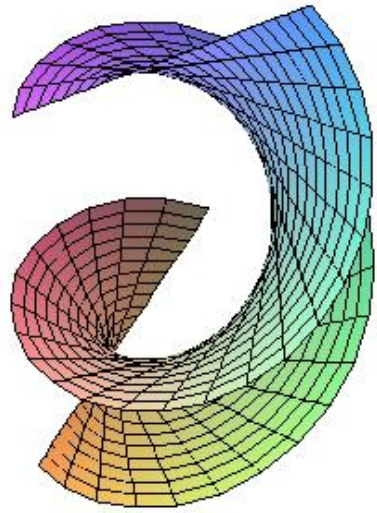
**Fig. 2.** Similar surface $N_\alpha$ of $N_\beta$.

## 5. Conclusions

A family of ruled surfaces in the Euclidean 3-space $E^3$ are defined and called similar ruled surfaces. Some properties of these special surfaces are obtained and it is showed that developable ruled surfaces form a family of similar ruled surfaces if and only if the striction curves of the surfaces are similar curves with variable transformation. By considering the importance of the offset surfaces we hope this paper leads new characterizations of ruled surfaces in different spaces.

**References**

[1] El-Sabbagh, M.F., Ali, A.T.: Similar Curves with Variable Transformations, Konuralp J. of Math., Vol. 1, No. 2, 2013, 80-90.

[2] Karger, A., Novak, J.: Space Kinematics and Lie Groups. STNL Publishers of Technical Lit., Prague, Czechoslovakia (1978).

[3] Küçük, A., Gürsoy O.: On the invariants of Bertrand trajectory surface offsets, App. Math. and Comp., 151 (2004) 763-773.

[4] Liu, H., Wang, F.: Mannheim partner curves in 3-space, Journal of Geometry, vol. 88, no. 1-2, pp. 120-126, 2008.

[5] Orbay, K., Kasap, E., Aydemir, İ.: Mannheim Offsets of Ruled Surfaces, Mathematical Problems in Engineering, Volume 2009, Article ID 160917.

[6] Papaioannou, S.G., Kiritsis, D.: An application of Bertrand curves and surfaces to CAD/CAM, Computer Aided Design, 17 (8) (1985) 348-352.

[7] Pottmann, H., Lü, W., Ravani, B.: Rational ruled surfaces and their offsets, Graphical Models and Image Processing, 58 (6) (1996) 544-552.

[8] Ravani, B., Ku, T.S.: Bertrand Offsets of ruled and developable surfaces, Comp. Aided Geom. Design, 23 (2) (1991) 145-152.

[9] Struik, D.J.: Lectures on Classical Differential Geometry, 2$^{nd}$ ed. Addison Wesley, Dover, (1988).

[10] Wang, F., Liu, H.: Mannheim partner curves in 3-Euclidean space, Mathematics in Practice and Theory, vol. 37, no. 1, pp. 141-143, 2007.